\documentclass[a4paper, 12pt]{article}
\usepackage[margin=2.5cm]{geometry}
\usepackage[utf8]{inputenc}
\usepackage[nosumlimits]{amsmath}
\usepackage{amssymb,stmaryrd,amsthm,MnSymbol}
\usepackage{amsfonts,amssymb,amsopn}
\usepackage[colorlinks=true,
            linkcolor=red,
            filecolor=blue,
            citecolor=blue,
            urlcolor=blue,]{hyperref}
\usepackage{color}
\usepackage[dvipsnames]{xcolor}
\usepackage{comment}
\usepackage{graphicx}
\usepackage{natbib}
\usepackage{multirow}

\usepackage{tikz}

\usepackage{soul}
\setstcolor{red}

\title{Fully implicit timestepping methods for the rotating shallow water equations}
\author{Werner Bauer\footnote{Corresponding author:
w.bauer@surrey.ac.uk,
School of Mathematics and Physics,
University of Surrey,
Guildford, GU2~7XH, UK,
ORCID:~\href{https://orcid.org/0000-0002-5040-4287}{https://orcid.org/0000-0002-5040-4287}
}
 ,
Colin J. Cotter\footnote{Imperial College London, UK, ORCID:
\href{https://orcid.org/0000-0001-7962-8324}{https://orcid.org/0000-0001-7962-8324}}
, and
Joshua Hope-Collins\footnote{Imperial College London, UK, ORCID:
\href{https://orcid.org/0000-0003-2051-8627}{https://orcid.org/0000-0003-2051-8627}}
}

\DeclareMathOperator{\diff}{d}
\DeclareMathOperator{\gat}{\diff\!}

\newcommand{\jump}[1]{[\![#1]\!]}

\begin{document}

\maketitle

\begin{abstract}
  
Fully implicit timestepping methods have had limited application in
geophysical fluid dynamics due to challenges of finding suitable
iterative solvers, since the coupled treatment of advection prevents
the standard elimination techniques. However, overlapping Additive
Schwarz methods provide a robust, scalable iterative approach for
solving the monolithic coupled system for all fields and Runge-Kutta
stages.
In this study we investigate this approach applied to the rotating
shallow water equations, facilitated by the Irksome package which
provides automated code generation for implicit Runge-Kutta
methods. We compare various schemes in terms of accuracy and
efficiency using an implicit/explicit splitting method, namely the
ARK2 scheme of Giraldo et al (2013), as a benchmark. {This provides an
initial look at whether implicit Runge-Kutta methods can be viable for
atmosphere and ocean simulation.}


\end{abstract}

\section{Introduction}

The selection of timestepping algorithms for numerical methods for
atmosphere models (and for large scale numerical solutions of partial
differential equations in general) is a complex issue due to the
interplay between the timestepping method, the timescales present, the
spatial discretisation, the types of implicit systems to be solved,
the types of solver algorithms that are available, the parallel
scalability and efficiency of those algorithms, \emph{etc.} For
atmosphere models, there are a wide range of approaches to
timestepping. Semi-implicit semi-Lagrangian methods are used by two of
the leading operational centres, the European Centre for Medium Range
Weather Forecasts \citep{hortal2002development} and the Met Office
\citep{melvin2010inherently,wood2014inherently}, with the next
generation Met Office model using a variant based on finite volume
methods for the advective transport
\citep{melvin2019mixed,melvin2024mixed}. Other approaches include
split explicit methods \citep{klemp2007conservative} such as those
used in the Weather Research and Forecasting (WRF) model
\citep{skamarock2008description} and the German Weather Service ICON
model \citep{zangl2015icon}, and implicit-explicit (IMEX) schemes
\citep{vogl2019evaluation,giraldo2013implicit}.

In this article we explore another direction which is fully implicit
Runge-Kutta methods. These methods avoid the complications of
splitting methods (stability conditions, options of composition, etc)
and bring in new complications of solving the resulting sparse
implicit systems. The main challenge is that the semi-implicit schemes
above involve the solution of linear implicit systems that can be
reduced down to a single scalar elliptic equation to be solved using
classical iterative strategies. For fully implicit schemes implemented
using Newton's method, the full Jacobian is involved; it cannot be
reduced down to a single variable due to the advective terms. Hence,
in this article we investigate monolithic additive Schwarz methods
that tackle the full system of variables. In the present work we make
an initial investigation of this strategy applied to the rotating
shallow water equations on the sphere. This system is two dimensional
so does not require a large facility to conduct experiments. Our
implementation, facilitated by the Firedrake
\citep{FiredrakeUserManual}, Irksome \citep{farrell2021irksome} and
PETsc \citep{balay2019petsc} software packages, uses a compatible
finite element spatial discretisation \citep{cotter2023compatible},
but this is not the central focus of the work.

The rest of the paper is organised as follows. In Section \ref{sec:methods},
we describe the space and time integration methods and the iterative
solver algorithm which we are using. In Section \ref{sec:results} we
present results of numerical experiments that compare the schemes for
speed and accuracy with the ARK2 IMEX scheme, a well known approach that
will allow the community to relate the performance to other methods.
In Section \ref{sec:summary}, we provide a summary and outlook.

\section{Description of methods}
\label{sec:methods}
\subsection{Spatial discretisation}
In this article we consider implicit time discretisations for the
rotating shallow water equations on the sphere, which we write
here in vector-invariant form as
\begin{align}
  u_t + \left(\nabla^\perp\cdot u + f\right)u^\perp
  + \nabla\left(\frac{|u|^2}{2} + g(D + b)\right) & = 0, \\
  D_t + \nabla\cdot(uD) & = 0,
\end{align}
where $u$ is the velocity (tangential to the sphere), $\nabla$ is the
gradient projected into the tangent plane on the sphere, $v^\perp =
k\times v$ for vector fields $v$, $\nabla^\perp\cdot u = k \cdot \nabla\times u$,
$k$ is the unit outward pointing
normal to the sphere, $f=2\Omega \sin(\phi)$ is the Coriolis parameter
with $\Omega$ the rotation rate of the Earth (2$\pi$/(sidereal day))
and latitude $\phi$, $g$ is the acceleration due to gravity, $b$ is
the topography field, and $D$ is the depth of the layer.

In this investigation we use a compatible finite element
discretisation of these equations. We do not expect that the precise
details are important for our conclusions, which are hopefully
  translatable to other discretisation approaches. However, to
efficiently describe our iterative solver approach a precise
description is useful. We select $V$ as the degree $p+1$ BDM
(Brezzi-Douglas-Marini) finite element space on triangles,
and $Q$ as the degree $p$ discontinuous Lagrange finite element space,
here defined on an icosahedral grid $\mathcal{X}$ formed by
recursively refining an icosahedron and then projecting vertices
radially out to the sphere.  Then we seek $(u,D)\in V\times Q$ such
that
\begin{align}
  \label{eq:ut}
  \langle w, u_t \rangle + a(u,D;w) 
   & = 0,
  \quad \forall w \in V, \\
  \label{eq:Dt}
  \langle \phi, D_t \rangle
+ c(u,D; \phi) & =
 0, \quad \forall \phi \in Q,
\end{align}
where
\begin{align}
   a(u,D;w) 
&=
  \langle w, fu^\perp \rangle 
  - \langle \nabla_h^\perp (w\cdot u^\perp), u \rangle
  + \llangle \jump{(w\cdot u^\perp) n^\perp}, \tilde{u} \rrangle  - \left\langle \nabla\cdot w, \frac{|u|^2}{2} + g(D
  + b) \right\rangle, \\
  c(u,D; \phi) &=
  - \langle \nabla_h \phi, uD \rangle
  + \llangle \jump{\phi n\cdot u}, \tilde{D} \rrangle,
\end{align}
and where $\langle \cdot , \cdot \rangle$ is the usual $L^2$ inner
product defined for scalar or vector fields integrating over the
domain $\mathcal{X}$, $\llangle\cdot,\cdot \rrangle$ is the $L^2$
inner product integrating over the set $\Gamma$ of mesh facets,
$\nabla_h$ is the ``broken'' cellwise gradient,
$\nabla_h^\perp \chi = -\nabla_h \times (k\chi)$ for a scalar field $\chi$,
$n$ is the unit normal (bivalued so that on each side of the facet $n$ is oriented to point
into the other side),
$\jump{\psi}$ indicates the sum of the
values of $\psi$ over both sides of the facet\footnote{To
  explain the $\jump{\cdot}$ operator further, we consider labelling
  the two sides of each edge with $+$ and $-$, so that $\psi^-$
  indicates the value of $\psi$ on the $-$ side, etc., thus
  $\jump{\psi}=\psi^++\psi^-$. Note that when $n$ is present,
  $n^+=-n^-$, and hence e.g. $\jump{\phi n \cdot u}=\phi^+n^+\cdot u^+
  + \phi^{-}n^-\cdot u^- = n^+(\phi^+u^+ - \phi^-u^-)$, so this gives
  rise to a jump. This particular choice of notation is used because
  of the appearance of $n^\perp$, which does not match standard
  patterns in mixed formulations of elliptic problems or discontinuous
  Galerkin discretisations for scalar transport.}, and
$\tilde{u}$ and
$\tilde{D}$ are the values of $D$ and $u$ evaluated on the upwind side
of a facet (the side with $u\cdot n > 0$). Here, the
  semicolons in the forms $a(u,D;w)$ and $c(u,D;\phi)$ indicate the
  linearity of the operators in the test functions $w$ and $\phi$,
  respectively. The forms can be nonlinear in the solution variables;
  this is a standard formulation for finite element methods applied
  to nonlinear partial differential equations.
For more details, and a
derivation of the finite element spaces and this finite element
approximation, see \cite{gibson2019compatible}. In this work we used
$p=1$.

\subsection{Implicit Runge-Kutta time discretisation and the resulting implicit systems}

Runge-Kutta methods advance the solution $U$ from one step $U^n$ to
the next $U^{n+1}$ by computing $s$ stages $(k_1,\ldots,k_s)$ and then
updating the solution according to
\begin{equation}
  U^{n+1} = U^n + \Delta t\sum_{i=1}^s b_i k_i,
\end{equation}
where $b=(b_1,\ldots, b_s)$ are coefficients specific to the
particular chosen Runge-Kutta method. Here, the stages are
  values of the time derivatives evaluated at different points in time
  relative to the current timestep; this is referred to as the ``stage
  derivatives'' formulation. There is an alternative ``stage values''
  formulation of Runge-Kutta methods where one solves for the value
  of the solution at different points instead, but we do not consider
  it here.

Our spatially discrete rotating
shallow water system (\ref{eq:ut}-\ref{eq:Dt}) is a ``mixed'' coupled
system for two variables, so we write $U=(u,D)\in V\times Q$ and $k_i
= (k_{u,i},k_{D,i})\in V\times Q$ for $i=1,\ldots,s$. Further,
implicit Runge-Kutta methods (IRKs) couple all of these stages together, so it
is useful to define a single variable for all of the stage components,
\[
k = (k_{u,1},k_{D,1},\ldots,k_{u,s},k_{D,s}) \in
\prod_{i=1}^s V\times Q := W.
\]
IRKs for (\ref{eq:ut}-\ref{eq:Dt}) then 
seek $k\in W$ such that
\begin{align}
    \label{eq:u stages}
    R_{u,i}[k;w] := \langle w, k_{u,i} \rangle +  a\left(
    u^n_i,D^n_i;w\right)
  & = 0,  
  \quad \forall w \in V,\mbox{ for }i=1,\ldots,s, \\
  \label{eq:D stages}
  R_{D,i}[k;\phi] := \langle \phi, k_{D,i} \rangle +  c\left(u^n_i,D^n_i; \phi\right)
  & = 0, \quad \forall \phi \in Q, \mbox{ for }i=1,\ldots,s,\\
  u^n_i = u^n + \Delta t\sum_{j=1}^sA_{ij}k_{u,j}, & \mbox{ for }i=1,\ldots,s,\\
  D^n_i = D^n + \Delta t\sum_{j=1}^sA_{ij}k_{D,j}, & \mbox{ for }i=1,\ldots,s,
\end{align}
where $A_{ij}$ are the matrix coefficients from the Butcher
tableau for the chosen Runge-Kutta scheme. This defines a coupled
system for all of the stages in general.

In this work
we consider collocation Runge-Kutta methods, in particular the
Gauss-Legendre methods that extend the implicit midpoint rule to
higher orders, and the Radau IIA methods that extend the backward Euler
method to higher orders. The Gauss-Legendre methods are energy
preserving for wave equations but slow down high frequency
oscillations, whilst the Radau IIA methods damp high frequency
oscillations.  See the two volume textbook \citet{hairer1993solving,wanner1996solving} for a comprehensive
derivation along with their Butcher
tableau (Chapter II.7), and stability and convergence analysis (Chapters IV.3-IV.5).
Note here that our system has no explicit time dependence
which would otherwise need to be incorporated into (\ref{eq:u
  stages}-\ref{eq:D stages}) in the usual way.

In this article we will discuss Newton-type methods for solving
 (\ref{eq:u stages}-\ref{eq:D
  stages}). Given an initial guess $k\in W$ for
the stages, a Newton iteration requires solving the coupled linear
Jacobian system for corrections
\[
k' = (k'_{u,1},k'_{D,1},\ldots, k'_{u,s},k'_{D,s}) \in W,
\]
such that
\begin{align}
  \nonumber
    \langle w, k'_{u,i} \rangle +
   \Delta t \sum_{j=1}^sA_{ij}\gat a\left(u^n_i,D^n_i;(k'_{u,j}, k'_{D,j}),
   w\right) &=  -R_{u,i}[k;w], \\
   \label{eq:u J}
   & \qquad\qquad    \forall w \in V,\mbox{ for }i=1,\ldots,s, \\
   \nonumber
  \langle \phi, k'_{D,i} \rangle
  + \Delta t\sum_{j=1}^sA_{ij}\gat c\left(u^n_i,D^n_i;(k'_{u,j}, k'_{D,j}),
  \phi\right)  &= -R_{D,i}[k;\phi], \\
  & \qquad \qquad \forall \phi \in Q,\mbox{ for }i=1,\ldots,s.
    \label{eq:D J}
\end{align}
Here, $\gat a$ and $\gat c$ are the Gateaux
derivatives of $a$ and $c$ defined by
\begin{align}
  \gat a(u,D;(v,\phi),w) & = \lim_{\epsilon\to 0}
  \frac{1}{\epsilon}(a(u+\epsilon v,D+\epsilon \phi; w)
  - a(u,D;w)),  \, \forall u,v,w \in V, \, D,\phi \in Q, \\
  \gat c(u,D;(v,\phi),\psi) & = \lim_{\epsilon\to 0}
  \frac{1}{\epsilon}(c(u+\epsilon v,D+\epsilon \phi; \psi)
  - c(u,D;\psi)), \, \forall u,v \in V,\, D,\phi,\psi \in Q.
\end{align}
In this work, we
take the convention that the upwind switches have derivative
zero\footnote{This treatment means that we {are already considering a
quasi-Newton method here}; using a semi-smooth Newton method might be
more elegant, but we do not observe any issues with our approach
when $u\cdot n=0$.}.
The solution is then updated according to
\begin{equation}
  k_{u,i}\mapsto k_{u,i} + k'_{u,i}, \quad
  k_{D,i}\mapsto k_{D,i} + k'_{D,i},
  \label{eq:newton update}
\end{equation}
and the iteration is repeated until the residuals in (\ref{eq:u
  stages}-\ref{eq:D stages}) are sufficiently small (according to an
appropriately chosen termination criteria; we use a relative
  tolerance test requiring that the ratio of the norm of the residual
  to the residual coming from the initial guess is below a selected
  tolerance, given later). The assembly of these nonlinear and linear
systems can be performed in the usual way by looping (in parallel)
over cells and constructing the appropriate contributions upon
substituting $w$ and $\phi$ for each basis function from standard
sparse finite element bases for $V\times Q$. It remains to find a
scalable way to solve (\ref{eq:u J}-\ref{eq:D J}).

\subsection{Scalable iterative methods}

  In this section, we describe our iterative approach to
  solving (\ref{eq:u stages}-\ref{eq:D stages}). This involves
  the composition of several stages, which we illustrate in Figure
  \ref{fig:solver_diag}.

\subsubsection{Nonlinear solver}
  
  At the outer level, (\ref{eq:u stages}-\ref{eq:D stages}) is solved
  using a quasi-Newton method, each iteration consisting of the
  approximate solution of the Jacobian system (\ref{eq:u J}-\ref{eq:D J}) followed by the
  update \eqref{eq:newton update}. These iterations are repeated until
  a relative tolerance condition is achieved, requiring that the norm
  of the residual in (\ref{eq:u stages}-\ref{eq:D stages}) is below
  some chosen tolerance multiplied by the residual norm before the
  first iteration. The choice of tolerance for our experiments is
  discussed in Section \ref{sec:results}.

  In the quasi-Newton method, an iterative method is used to solve
  (\ref{eq:u J}-\ref{eq:D J}), which is completed when a termination
  criteria is satisfied.  
  Our quasi-Newton method uses the Eisenstat-Walker (version 2)
  inexact Newton approach \citep{eisenstat1996choosing} which
  adaptively controls the {relative tolerance for iterative
    solution of Jacobian system (\ref{eq:u J}-\ref{eq:D J}) at each
    Newton step, since it is not useful to solve the Jacobian
  system accurately when the nonlinear solution is only going to be
  updated again in the next Newton iteration anyway. Here, relative
  tolerance means a maximum value of the ratio of the equation residual
  divided by the value before the first iteration, below which the
  Jacobian solution procedure is terminated.
  To avoid over solving the Jacobian system during the final Newton
  iteration, we also set an absolute tolerance equal to $0.1\times
  \tau$ multiplied by the norm of the starting nonlinear residual,
  where $\tau$ is the chosen nonlinear tolerance. If the true
  (\emph{i.e.} unpreconditioned) residual falls below this tolerance,
  then the Jacobian system iteration
  is terminated.  This works because
  near convergence of the nonlinear solver, the true residual for the
  Jacobian system is very close to the true residual for the nonlinear
  system. Hence, if the true residual for the Jacobian system is 10
  times less than the tolerance set for the nonlinear residual, then
  we expect that the nonlinear tolerance termination criteria to be
  satisfied on the start of the next iteration.

  Finally, we also ``lagged'' the Jacobian, so that it is only updated
  every 20 Newton iterations (including between Newton solves). This
  reduces the cost of rebuilding factorisations in the patch
  smoothers, at the expense of potentially increasing the number of
  required Newton iterations.

\subsubsection{Krylov method for Jacobian system}

We observe that the stage components of the iterative corrections $k'$
are coupled in (\ref{eq:u J}-\ref{eq:D J}) both between $u$ and $D$
components for the same stage $i$ (where the coupling is global but
sparse), as well as between all of the stages from $i=1$ to $s$. This
presents challenges for the solver. In this work, we solve (\ref{eq:u
  J}-\ref{eq:D J}) using a monolithic approach, meaning that we apply
a preconditioned Krylov solver to the full set of basis coefficients
for the iterative corrections $k'$, treated as a single vector. In our
experiments we used FGMRES (the flexible generalised minimum residual
method \citep{saad1993flexible}) for the chosen Krylov
method. This is because we used GMRES within the preconditioner (see later), which necessitates a flexible
  Krylov method outside.}

\begin{figure}
  \center
  \includegraphics[width=5cm]{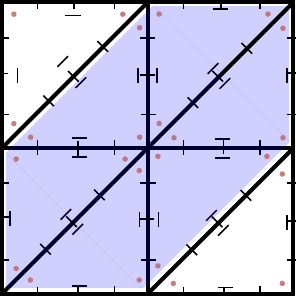}
  \caption{\label{fig:patch} A diagram showing the ``vertex star''
    patch associated with the vertex at the middle of the
    diagram. Dots indicate $D$ degrees of freedom and thin lines
    indicate $u$ degrees of freedom. Of the latter, lines crossing
    cell edges indicate normal components of $u$, which are continuous
    across cell edges for our choice of finite element space.  Lines
    parallel to cell edges indicate tangential components of $u$,
    which are discontinuous across cell edges.  The shaded region
    indicates the patch.  All degrees of freedom outside the shaded
    region (including normal but not tangential components on the
    patch boundary) are excluded from the patch. This
      diagram shows the degrees of freedom for a single stage, so for
      an $s$-stage scheme the number of degrees of freedom in the
      implicit coupled system is $s$ times bigger.  }
 \end{figure}

\subsubsection{Preconditioner for the Krylov method}

As a preconditioner, we used geometric multigrid with a sequence of
nested meshes of the spherical domain,\footnote{In fact our meshes are
not strictly nested, because our cells are flat triangles (higher
order polynomial representations are also possible). Starting from the
coarsest icosahedral mesh, each triangle is refined by replacing it
with four smaller triangles, created by adding vertices at the
midpoint of each edge. Hence, the meshes can be described as
  being topologically nested but not geometrically nested.
The new vertices are then moved to the surface
of the sphere being approximated. Transfer operators are defined by
moving the extra vertices of the finer mesh to the edges of the coarse
mesh before prolongation/restriction.} with prolongation operators
consisting of the inclusion operator, \emph{i.e.} reinterpretation of
the same function in the larger space on the finer mesh; restriction
is taken as the $L_2$ dual of prolongation.

For smoothers on each level above the coarsest,
we used an additive Schwarz method
(ASM) that
is built using overlapping subspaces $W_l\subset W
{=\prod_{i=1}^s\left(V \times Q\right)}$, $l=1,\ldots,N_{{L}}$,
where $W$ is the mixed finite element space comprising all $s$ stages on multigrid level $L$,
  and where $N_L$ is the number of overlapping subspaces.
In this work we define the subspaces $W_l$ using ``vertex
star patches''.
Here, there is one subspace $W_l$ for each of the $M$
vertices $z_l$ in the mesh.
The subspace for a vertex $z_l$ is defined by
taking the set $S_l$ of cells surrounding $z_l$. $W_l$ is the subspace
of $W$ consisting of functions that are equal to zero when restricted
to any cell not in $S_l$. For our choice of spaces, this entails
zeroing any degrees of freedom on the boundary of $S_l$ and beyond
(see Figure \ref{fig:patch}). In the multigrid algorithm, the goal of
the smoother is to approximately solve a problem of the form
(\ref{eq:u J}-\ref{eq:D J}) but with $R_{u,i}$ and $R_{D,i}$ replaced
by appropriate linear residuals.

For the additive
Schwarz method, we solve (\ref{eq:u J}-\ref{eq:D J}) with solutions
$w_l$ and test functions restricted to $W_l$, using a direct
solve. These ``patch problems'' can be solved independently and in
parallel using a dense direct solver.  Then, the additive Schwarz
approximate solutions are obtained as $\sum_{l=1}^{N_{{L}}}w_l$ (here
$w_l$ is interpreted as a function in $W$ by inclusion $W_l\subset
W$).
\begin{figure}
\tikzstyle{decision} = [rectangle, draw, fill=blue!20, 
    text width=4.5em, text badly centered, node distance=3cm, inner sep=0pt]
\tikzstyle{block} = [rectangle, draw, fill=blue!20, 
    text width=5em, text centered, rounded corners, minimum height=4em]
\tikzstyle{line} = [draw, -latex]

\scalebox{0.7}{
\begin{tikzpicture}[node distance = 3cm]
  \node [block] (NL) {Nonlinear system};
    \node [block, below of=NL] (Newton) {Eisenstat Walker Quasi-Newton};
    \node [block, right of=Newton, node distance=6.5cm] (FGMRES) {FGMRES};
    \node [block, right of=FGMRES, node distance=3.5cm] (MG) {Multigrid};
    \node [block, right of=MG, node distance=6.8cm] (Coarse) {GMRES 3 its};
    \node [block, right of=Coarse, node distance=3.5cm] (CPC) {Block ILU};
    \node [block, below of=Coarse] (Levels) {GMRES 3 its};
    \node [block, right of=Levels, node distance=3.5cm] (LPC) {Patch ASM};
    \path [line] (NL) -- node[right]{Nonlinear solver}(Newton);
    \path [line] (Newton) -- node[below]{linear Jacobian solver}(FGMRES);
    \path [line] (MG) -- node [below]{Coarse grid smoother} (Coarse);
    \path [line] (Coarse) -- node [below]{PC} (CPC);
    \path [line] (Levels) -- node [below]{PC} (LPC);
    \path [line] (MG) |-node [below right]{Levels smoother}  (Levels);
    \path [line] (FGMRES) -- node [below]{PC} (MG);
\end{tikzpicture}
} \caption{\label{fig:solver_diag}Diagram showing the
    composition of solver components for our approach. The boxes
    indicate solvers which have components indicated by the
    arrows. Here PC indicates a preconditioner.}
\end{figure}

The reason for choosing vertex star patches is that they deal well
with the oscillatory wave coupling in the linear system. When this
form of multigrid is applied to the rotating shallow water equations
linearised about a state of rest (referred to as the linear rotating
shallow water equations), we observe mesh- and $\Delta t$-robust
convergence rates.
This is believed to be related to the efficacy of
vertex star patches for Hdiv problems \citep{arnold2000multigrid} but
there is no analysis for monolithic multigrid applied to mixed
elliptic problems at present \citep{sm}. The reason for choosing patches
that couple between all stages is that we aim to obtain robust
convergence in the number of Runge-Kutta stages, as analysed in
\citep{kirby2024convergence}.

To avoid having to tune scaling parameters whose optimal values might
depend on the system state, our smoothers consist of 3
iterations of GMRES (the generalised minimum residual
  method \citep{saad1986gmres}) preconditioned by the additive
Schwarz method above. This has the effect of selecting the scaling
parameter adaptively.  This necessitates the use of a flexible Krylov
method (FGMRES, in our case) for the ``outer'' monolithic solver,
since the use of GMRES here means that the smoother is residual
dependent and hence is not a stationary iterative method (which is a
requirement for standard Krylov methods such as GMRES).

In our multigrid setup, we used 3 iterations of GMRES
  preconditioned by block ILU (incomplete LU factorisation) for the
correction on the coarsest grid. This avoids having to use a parallel
direct solve. We experimented with a direct solve on the coarse grid
but found that it did not alter the overall convergence of the solver
strategy, but it did increase wallclock time.

As will become apparent from the results, we do not observe perfect
multigrid behaviour in that the number of iterations is not robust in
$\Delta t$: more iterations are required at larger $\Delta t$ due to
the presence of the advective terms.  For this type of problem, it is
also typical to scale $\Delta t$ to keep the advective Courant number
($U\Delta t/\Delta x$ where $U$ is a typical velocity scale and
$\Delta x$ is a typical cell diameter) constant as the mesh is
refined, so we can hope for mesh independent Krylov iteration counts
under this refinement scaling without multigrid (as obtained
in \citet{cottershipton2023compatible}).
However, we found
that using the multigrid produced faster wallclock times than
simply applying the patch smoother on the finest grid.

\subsubsection{Solver software implementation}

  In this work we implemented this approach using Irksome
  \citep{farrell2021irksome,kirby2024extending}, a Python library that
  wraps time discretisations on top of finite element spatial
  discretisations implemented using Firedrake
  \citep{FiredrakeUserManual}, an automated system for the solution of
  partial differential equations using the finite element method.
  This combination in turn makes significant use of PETSc
  \citep{dalcin2011parallel,balay2019petsc}. In particular, PETSc's
  and Firedrake's PCPatch implementation of additive Schwarz methods
  is used \citep{farrell2021pcpatch}, via the Python preconditioner
  interface of \citet{kirby2018solver}.  Amongst the benefits of this
  approach to implementation, Firedrake uses automated differentiation
  provided by the Unified Form Language (UFL) \citep{alnaes2012ufl} to
  obtain the formulae for (\ref{eq:u J}-\ref{eq:D J}) from which code
  assembling the matrix-vector action and the patch linear systems is
  automatically generated. In our implementation, we use $\ell_2$
  Euclidean norms of vectors of finite element basis coefficients to
  measure residual norms for termination criteria.

\subsection{Implicit-explicit time discretisation}

We compared our IRK approach with an
implicit-explicit (IMEX) time discretisation. This serves as a
benchmark in terms of timing since IMEX methods are well known in the
community developing numerical methods for geophysical fluid dynamics.
We choose one specific second-order IMEX scheme, the ARK2 scheme
of \citet{giraldo2013implicit}, which they demonstrated to be optimal
for geophysical fluid problems amongst second order IMEX schemes and has
been widely adopted in the field.

To produce an IMEX scheme, we rewrite our spatial discretisation
(\ref{eq:ut}-\ref{eq:Dt}) in the form
\begin{align}
  \label{eq:ut IMEX}
  \langle w, u_t \rangle + a_L(u,D;w) 
  + a_N(u,D;w) 
   & = 0,
  \quad \forall w \in V, \\
  \label{eq:Dt IMEX}
  \langle \phi, D_t \rangle
+ c_L(u,D; \phi) + c_N(u,D; \phi)  & =
 0, \quad \forall \phi \in Q,
\end{align}
where $a_L$, $a_N$ are the linear and nonlinear operators for $u$,
and $c_L$, $c_N$ are the linear and nonlinear operators for $D$,
defined by
\begin{align}
  a_L(u,D;w) & = \langle w, fu^\perp \rangle - \langle \nabla\cdot w,
  gD \rangle, \\
   a_N(u,D;w) 
&=
  - \langle \nabla_h^\perp (w\cdot u^\perp), u \rangle
  + \llangle \jump{(w\cdot u^\perp) n^\perp}, \tilde{u} \rrangle  - \left\langle \nabla\cdot w, \frac{|u|^2}{2} + gb \right\rangle, \\
  c_L(u,D;\phi) & = \langle \phi, H\nabla\cdot u\rangle, \\
  {c_N}(u,D; \phi) &=
  - \langle \nabla_h \phi, uD \rangle
  - \langle \phi, H\nabla\cdot u\rangle
  + \llangle \jump{\phi n\cdot u}, \tilde{D} \rrangle,
\end{align}
where $H$ is the depth $D$ at rest. Here, the linear terms $a_L$ and
$c_L$ are from the linearisation of the equations at a state of rest,
which have fast wave solutions.  The nonlinear remainder terms $a_N$
and $c_N$ are advective and contribute on the slower advective
timescale. Hence, it can be beneficial to use this splitting in an
IMEX scheme where the linear terms are treated implicitly.

In our setting, IMEX schemes take the form,
\begin{align}
    \nonumber
    \langle w, Y_{u,i}-u^n \rangle + \Delta t \sum_{j=1}^{i-1} A_{ij} a_N\left(
    Y_{u,j}, Y_{D,j}; w\right) & \\
\qquad    + \Delta t \sum_{j=1}^i \tilde{A}_{ij} a_L(Y_{u,j}, Y_{D,j}; w)
    & = 0, \label{eq:u stages IMEX}
  \quad \forall w \in V,\mbox{ for }i=1,\ldots,s, \\
  \nonumber
    \langle \phi, Y_{D,i}-D^n \rangle + \Delta t \sum_{j=1}^{i-1} A_{ij} c_N\left(
    Y_{u,j}, Y_{D,j}; \phi\right) & \\
    + \Delta t \sum_{j=1}^i \tilde{A}_{ij} c_L(Y_{u,j}, Y_{D,j}; {\phi})
    & = 0, \label{eq:D stages IMEX}
      \quad \forall \phi \in Q,\mbox{ for }i=1,\ldots,s,
\end{align}
for stages $(Y_{u,i},Y_{D,i}) \in V \times Q$,
where $A$ and $\tilde{A}$ are the explicit and implicit IMEX Butcher
matrices respectively. Then, we reconstruct the solution from
\begin{align}
  \label{eq:mass u}
  \langle w, u^{n+1} - u^n\rangle + \Delta t\sum_{i=1}^s b_i
  a_N(Y_{u,i}, Y_{D,i}; w)
  + \Delta t \sum_{i=1}^s \tilde{b}_i
  a_L(Y_{u,i}, Y_{D,i}; w)=0,\, \forall w \in V, \\
  \label{eq:mass D}
    \langle \phi, D^{n+1} - D^n\rangle + \Delta t\sum_{i=1}^s b_i
  c_N(Y_{u,i}, Y_{D,i}; \phi)
  + \Delta t \sum_{i=1}^s \tilde{b}_i
    c_L(Y_{u,i}, Y_{D,i}; \phi)=0,\, \forall \phi \in Q,
\end{align}
where $b$ and $\tilde{b}$ are the corresponding IMEX Butcher
reconstruction vectors. Due to the lower triangular structure
expressed in the sum limits, (\ref{eq:u stages IMEX}-\ref{eq:D stages
  IMEX}) can be solved as single coupled mixed problems for
$(Y_{u,i},Y_{D,i})$ at stage $i$. Further, because the
  linear component of the splitting contains only the inertia-gravity
  wave system, these problems can be reduced to a sparse ``modified
Helmholtz'' type problem for a single variable using the hybridisation
technique for mixed finite elements
\citep{boffi2013mixed,cockburn2004characterization}; similar sparse
reduction techniques exist for finite volume and discontinuous
Galerkin methods. The resulting sparse system can be solved using a
sparse parallel direct solver, or for larger problems, multigrid
methods. In this work we found that sparse direct solvers were quicker
for the 2D problems under consideration. The implicit IMEX linear
systems are state independent, and so can be assembled and factorised
once, which makes their solution very fast. In contrast, the patch
problems used for the IRK methods depend on the system state, and must
be refactorised each time the nonlinear system Jacobian is
  recomputed. The principle disadvantage with the IMEX scheme is that
it is stable conditional on the advective Courant number, which means
that it is limited to smaller timesteps.  It is this property that we
shall contrast with the IRKs.

The ARK2 IMEX scheme is a 3 stage scheme involving two implicit solves
per timestep. It was implemented using Firedrake again, using the
hybridization solver package SLATE \citep{gibson2020slate} for the
implicit linear systems. The Butcher tableau for ARK2 are shown in Table
\ref{tab:ARK2}.

\begin{table}\centering
\begin{tabular}{c|ccc}
0         &  0             &   0       & 0  \\
$2\gamma$ &  $2\gamma$     &  0        & 0  \\
$1$       & $1 - \alpha$   & $\alpha$  & 0  \\ \hline
          & $\delta$       & $\delta$  & $\gamma$
\end{tabular}
\hspace{4mm}
\begin{tabular}{c|ccc}
0         &  0             &   0       & 0  \\
$2\gamma$ &  $\gamma$     &  $\gamma$        & 0  \\
$1$       & $\delta $   & $\delta$  & $\gamma$  \\ \hline
          & $\delta$       & $\delta$  & $\gamma$
\end{tabular}
\caption{\label{tab:ARK2}Double Butcher tableau for the ARK2 IMEX scheme of
  \citet{giraldo2013implicit}. Left: the tableau for the explicit component.
  Right: the tableau for the implicit component.
   Here, $\gamma = 1 -
  \frac{1}{\sqrt{2}}$, $\alpha = \frac{1}{6} (3 + 2\sqrt{2})$ and
  $\delta = \frac{1}{ 2 \sqrt{2}}$.}
\end{table}

\section{Numerical results}
\label{sec:results}

In this section, we present numerical experiments to investigate the
monolithic solver setup applied to IRKs for the rotating shallow water
equations as introduced above. Specifically, we will consider Gauss-Legendre and Radau IIA timestepping methods in their versions with 1,
2 and 3 stages, combined with the finite element spaces discussed
above in the case of $p=1$.
Additionally, we used the same solver approach with the 3 stage diagonally
   implicit Runge-Kutta method (DIRK) of \citet{alexander1977diagonally}.
We compare these methods with the
frequently used, very efficient IMEX method ARK2.
The goal is not to compete with the IMEX scheme, but rather to provide a
  point of reference using an approach that is familiar to the community,
  to give an idea of the viability of IRK schemes.

\begin{figure}
  \centering
 \begin{tabular}{cc}
 \hspace{-2em}\includegraphics[scale=0.3]{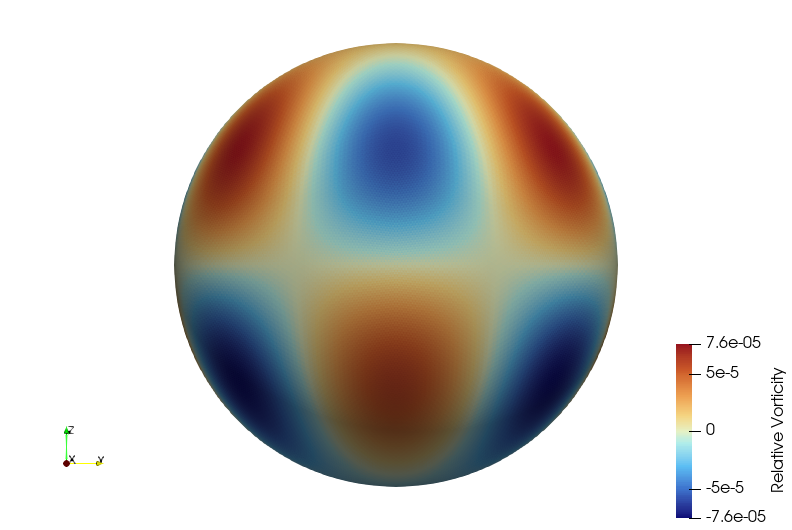} &
 \hspace{-1.5em}\includegraphics[scale=0.3]{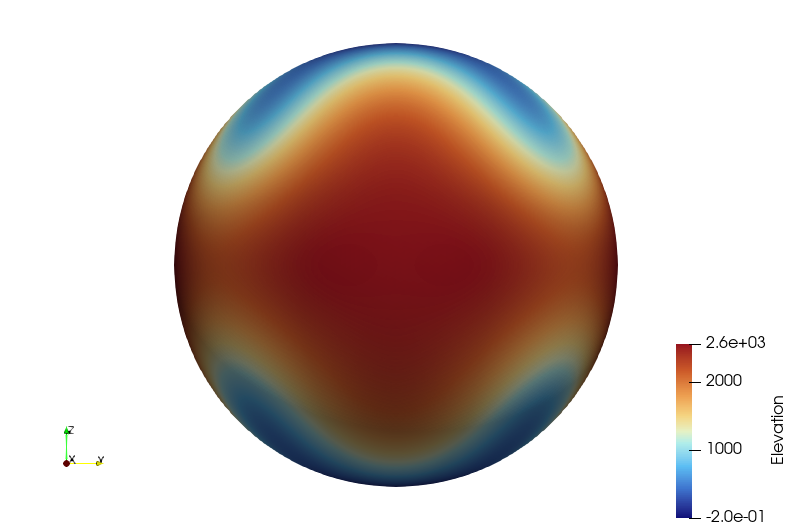} \\
 \hspace{-2em}\includegraphics[scale=0.3]{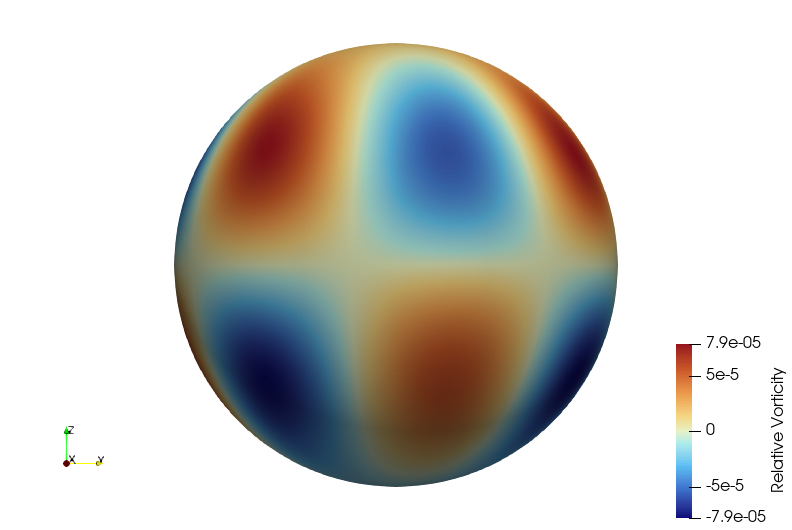} &
 \hspace{-1.5em}\includegraphics[scale=0.3]{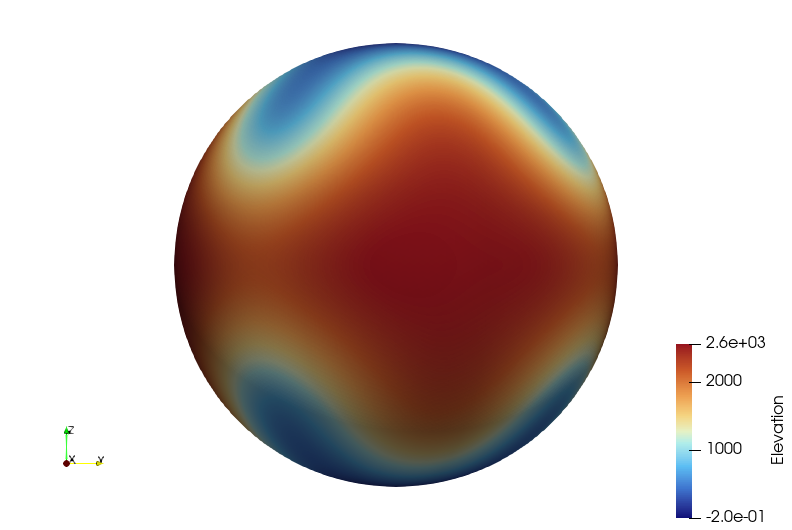} \\
 \end{tabular}\vspace{-10pt}
  \caption{Top row: initial vorticity (left column) and depth (right
    column) fields for Williamson test case 6. Bottom row: both fields
    at day 1. The shown fields are for mesh level 5, which has 20480
    cells, 61440 height degrees of freedom and 153600 velocity degrees
    of freedom. The reduced system in the IMEX scheme has 92160
    degrees of freedom.}
 \label{tc6_fields}
\end{figure}

The aim is to test this for geophysical applications, i.e for
atmosphere and ocean simulations, in the balanced regime where the
solution is dominated by the slow evolving component without
significant fast oscillations. To this end, we use the Rossby-Haurwitz
test case of \cite[Section 6]{WILLIAMSON1992211} which computes a
nonlinear travelling Rossby wave solution that is supported over the
whole globe. Plots of the solution are shown in Figure
\ref{tc6_fields}. We compute timestepping errors at relatively short time of 1 day.
The errors are computed relative to the IMEX numerical solution using the
same spatial discretisation scheme and mesh resolution with timestep
$\Delta t = 1s$. This enables us to isolate timestepping errors. We
also estimate the spatial discretisation error to give context to
these timestepping errors.

\begin{figure}\centering
 \begin{tabular}{c}
  \hspace{-1.3cm}\includegraphics[scale=0.62]{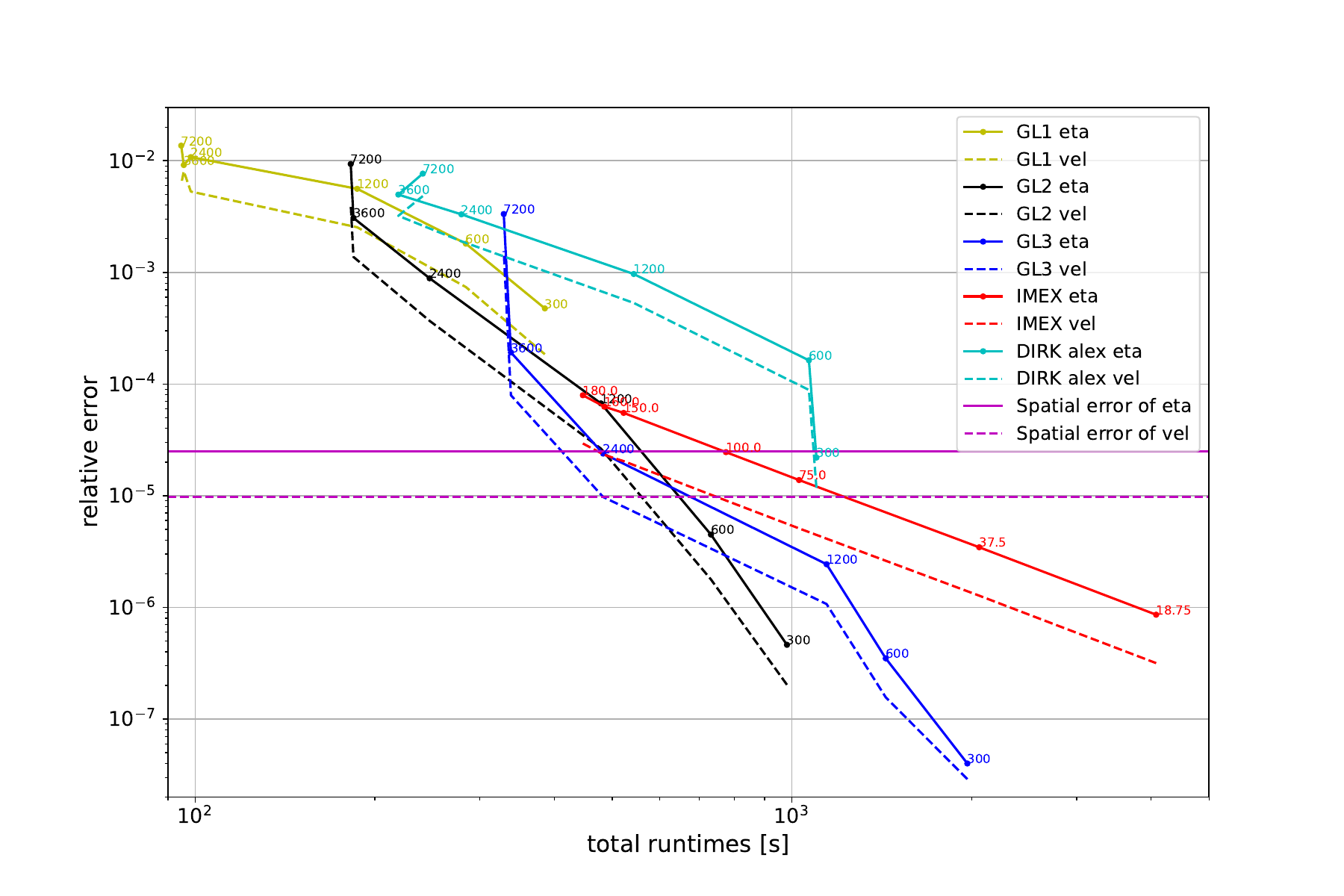}
 \end{tabular}
  \caption{Relative errors of Gauss-Legendre 1-3, Alexander DIRK, and ARK2 IMEX methods
    vs total wallclock time. The numbers on the curves indicate the
    time step sizes $\Delta t$ used for the corresponding solution. Solid
    lines indicate values for the elevation field $\eta=D-H$ and dashed lines
    for the velocity fields $ u$.
   The solid and dashed horizontal lines indicate an estimate of the spatial discretisation error between mesh levels L6 and L7
   of $\eta$ and $u$, respectively.
   }
   \label{err_vs_runtime_gl_imex}
\end{figure}

\begin{figure}\centering
\begin{tabular}{c}
 \hspace{-1.3cm}\includegraphics[scale=0.62]{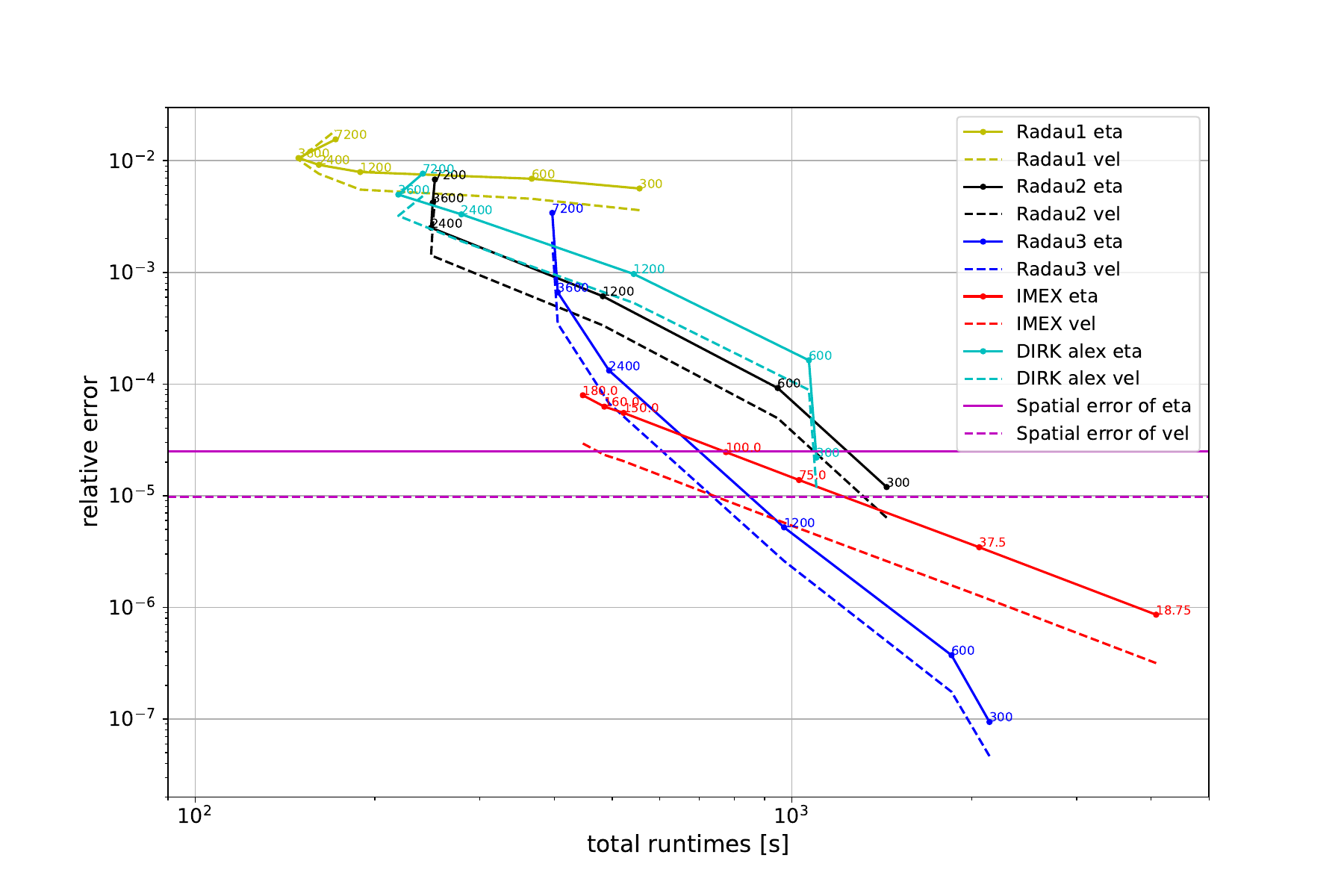}
\end{tabular}
   \caption{Relative errors of Radau IIA 1-3, Alexander DIRK, and ARK2 IMEX methods vs total
   wallclock time. The numbers on the curves indicate the
   time step sizes {$\Delta t$} used for the corresponding solution. Solid
   lines indicate values for the elevation field $\eta=D-H$ and dashed lines
   for the velocity fields $ u$.
   The solid and dashed horizontal lines indicate an estimate of the spatial discretisation error between mesh levels L6 and L7
   of $\eta$ and $u$, respectively.
   }
  \label{err_vs_runtime_rd_imex}
\end{figure}

As a benchmark to contextualise the timestepping errors, we
provide an indication of the spatial discretization error. This is
done by comparing the solution between mesh refinement levels L6 and
L7 after a 1 day simulation using a very small time step size of $\Delta t
=1s$ so that for the chosen mesh levels and time step size, the
temporal error is negligible compared to the spatial discretization
error.  We found spatial discretization errors of $2.542863 \cdot
10^{-05}$ for the depth field $D$ and $ 9.777218 \cdot 10^{-06}$
for the velocity field $u$. These error values are indicted as horizontal
lines in Figures \ref{err_vs_runtime_gl_imex} and
\ref{err_vs_runtime_rd_imex}. The point is that reducing the timestepping
error much below these values will not lead to an overall reduction in
the combined spatial and temporal error, so there is no benefit in
reducing the time stepping error further.

In our experiments, essential quantities of interest are the total runtime of a
simulation and the time errors the integrators have, up to the point
where they are of the same order of magnitude as the spatial
errors. As such, we compare these quantities for various parameter
choices for the IRK, DIRK, and IMEX schemes. As expected, accuracy and
runtime crucially depend on the selected time step sizes; we compute
errors and runtime for IRKs, DIRKs, and the ARK2 IMEX scheme over a range of
$\Delta t$ values. For the IMEX scheme, the maximum timestep was
$\Delta t=180s$ on mesh {refinement} 6 since the method was unstable
for larger timesteps.

All computations were performed using an icosahedral mesh
approximating the surface of the sphere by recursive refinement of an
icosahedron as discussed previously. {Mesh refinement level $p$}
 has $20\times 4^p$ cells, so there are $60\times 4^p$ layer depth degrees
 of freedom and $150\times 4^p$ velocity degrees of freedom. The
 machine used has 16 Intel Xeon 2.60GHz CPUs on two sockets, and we
 used one thread per core. All 16 cores were used unless stated
 otherwise, with timings on fewer cores presented in Table
 \ref{tab:cores}.  The timing measurements were obtained using PETSc's
 event logging capability.

\begin{table}[h]\centering
  \begin{tabular}{ccccccc}
    $\Delta t$ $(s)$ & \multicolumn{6}{c}{Wallclock time ($s$)} \\
    \hline
    &  16 cores  &  8 cores  & 4 cores  &   8/16   &   4/8     & 4/16     \\\hline
        \multicolumn{7}{c}{IMEX} \\
   180.00 &  437.0420 &  702.3965 &  1134.2820 & 1.60716 & 1.614874 & 2.595362\\
   18.75 &   4094.4577 & 6703.4870 & 10742.4779 & 1.63721 & 1.602521 & 2.623663 \\
    \multicolumn{7}{c}{Gauss-Legendre 1 stage} \\
  7200 &  91.961 & 134.849 & 224.160 & 1.466372 & 1.662304 & 2.437555\\
  1200 & 186.878 & 279.982 & 472.217 & 1.498207 & 1.686598 & 2.526873\\
   300 & 378.757 & 573.468 & 958.353 & 1.514079 & 1.671153 & 2.530258\\
  \multicolumn{7}{c}{Gauss-Legendre 2 stage} \\
  7200 & 181.891 &  270.150 &  456.909 & 1.485230 & 1.691316 & 2.511993\\
  1200 & 482.581 &  738.246 & 1247.484 & 1.529787 & 1.689794 & 2.585025\\
   300 & 979.752 & 1514.306 & 2567.877 & 1.545601 & 1.695745 & 2.620946\\
    \multicolumn{7}{c}{Gauss-Legendre 3 stage} \\
  7200 &  367.896 &  485.772 &  823.797 & 1.320406 & 1.695851 & 2.239212\\
  1200 & 1168.521 & 1790.578 & 3068.695 & 1.532346 & 1.713801 & 2.626136\\
   300 & 1957.424 & 3070.712 & 5300.002 & 1.568752 & 1.725985 & 2.707641\\
  \end{tabular}
  \caption{
  Table showing wallclock times for different numbers of cores for the ARK2 IMEX scheme and the Gauss-Legendre IRKs. Similar values are obtained for the Radau IIA IRKs. The results are typical for intra node parallelism, with the usual bandwidth limitations.
    \label{tab:cores}
    }
\end{table}

 For the IRKs we used the solver strategy discussed above with a
 Newton relative residual tolerance of $1.0\times 10^{-6}$,
 and a multigrid scheme with $p$ levels, where
   $p$ is the mesh refinement used.  For the
 IMEX scheme we used the SUPERLU Dist parallel direct solver package
 \citep{li2023newly,li2003superlu_dist} for the reduced hybridised
 system in the implicit solve since this was the fastest of those
 available through PETSc. The mass matrix system (\ref{eq:mass u}) was
 solved using GMRES preconditioned by ILU(0) (incomplete LU
 factorisation discarding all values outside the sparsity pattern of
 the mass matrix) to a tolerance of $1.0\times 10^{-10}$,
   and the mass matrix system (\ref{eq:mass D}) was solved using a
   direct solve since it is block diagonal. These systems have minor
   effect on the solve times compared to the two implicit solves per
 timestep.

   We also took measurements for the 3 stage DIRK scheme of
   \citet{alexander1977diagonally}, a third order method.  The
   diagonal structure means that the stages can be solved one by one
   instead of simultaneously, although a separate nonlinear system
   is still required to be solved at each stage.
   The system for each stage is solved using
   the same iterative solver framework as for the IRKs.
 
 Figures \ref{err_vs_runtime_gl_imex} and \ref{err_vs_runtime_rd_imex}
 show wallclock time versus error plots for the Gauss-Legendre and
 Radau IIA IRKs, and the Alexander DIRK; results for the ARK2 IMEX
 scheme are shown in both figures as a reference. Mesh refinement 6
 was used. We considered 1, 2 and 3 stage versions of the IRKs, as we
 were limited by memory on the workstation from examining more stages
 at this resolution. The Gauss-Legendre $s$-stage schemes are order
 $2s$ with stage order $s$ (the relevant order for stiff problems at
 large $\Delta t$), whilst the Radau IIA schemes are order $2s-1$ with
 stage order $s$.

 We observe that the IRKs and DIRKs allow for larger timesteps and
 hence quicker time to solution. The DIRK solution provides equivalent
 error at slower time than the IRKs: the collocation schemes provide
 higher accuracy for the same timestep, which is best seen by
 comparing the 2 stage Radau IIA scheme with the DIRK (they are both
 third order methods). We also see that the Gauss-Legendre schemes
 provide quicker time-to-solution for the same error than the Radau
 IIA schemes.  This is because the classical order of
   Gauss-Legendre schemes is one higher than that of Radau IIA schemes
   for the same number of stages.
 The
 plots also show that the IRKs provide a viable option since they are
 certainly not slower than the IMEX scheme in general.
 
 The loss of $\Delta t$ robustness for larger $\Delta t$ is demonstrated
 in Tables \ref{tab_niter_vs_dt G} and \ref{tab1_niter_vs_dt R}. We observe
 $\Delta t$ robustness for smaller $\Delta t$, with more sustained $\Delta t$
 robustness for the higher order schemes.

  We also studied the effect of changing the nonlinear solver
   tolerance from $10^{-6}$ to $10^{-5}$ and $10^{-7}$. There was
   insignificant changes to the relative solution error values
   (overall no bigger than $5\times 10^{-8}$), but there was a change in
   wallclock times, due to changes in the number of Newton
   iterations. When the tolerance was decreased, more Newton iterations
   were required in general, increasing the wallclock time (and similarly
   for increasing the tolerance).
   These changes are presented in Table \ref{tab:GLntol}.

\cleardoublepage

\begin{table}
  \centering
  \begin{tabular}{ccccccc}
    & \multicolumn{2}{c}{Gauss-Legendre 1} &
    \multicolumn{2}{c}{Gauss-Legendre 2} &
    \multicolumn{2}{c}{Gauss-Legendre 3} \\
   $\Delta t$  &  its/steps  &  its 
   &  its/steps  &  its
   &  its/steps  &  its 
    \\ \hline
7200 & 11.83 & 142.00 & 8.92 & 107.00 & 7.92 & 95.00 \\
3600 & 5.96 & 143.00 & 4.38 & 105.00 & 4.17 & 100.00 \\
2400 & 4.03 & 145.00 & 4.00 & 144.00 & 4.00 & 144.00 \\
1200 & 4.00 & 288.00 & 4.00 & 288.00 & 4.96 & 357.00 \\
600 & 3.00 & 432.00 & 3.01 & 433.00 & 3.01 & 433.00 \\
300 & 2.00 & 577.00 & 2.00 & 576.00 & 2.05 & 589.00 \\
\end{tabular}
\caption{Krylov iterations per step (its/steps) and total Krylov iteration count (its) vs $\Delta t$ for the Gauss-Legendre  1, 2, 3 methods.}
\label{tab_niter_vs_dt G}
\end{table}

\begin{table}[h]
  \centering
    \begin{tabular}{ccccccc}
        & \multicolumn{2}{c}{Radau IIA 1} &
    \multicolumn{2}{c}{Radau IIA 2} &
    \multicolumn{2}{c}{Radau IIA 3} \\
   $\Delta t$  &  its/steps  &  its 
   &  its/steps  &  its
   &  its/steps  &  its 
    \\ \hline
7200 & 23.08 & 277.00 & 12.83 & 154.00 & 9.83 & 118.00 \\
3600 & 9.71 & 233.00 & 6.21 & 149.00 & 5.04 & 121.00 \\
2400 & 6.94 & 250.00 & 4.06 & 146.00 & 4.11 & 148.00 \\
1200 & 4.03 & 290.00 & 4.00 & 288.00 & 4.12 & 297.00 \\
600 & 4.00 & 576.00 & 4.00 & 576.00 & 4.00 & 576.00 \\
300 & 3.00 & 864.00 & 3.00 & 864.00 & 2.25 & 647.00 \\
\end{tabular}
\caption{Krylov iterations per step (its/steps) and total Krylov iteration count (its) vs $\Delta t$ for the RADAU IIA 1, 2, 3 methods.}
\label{tab1_niter_vs_dt R}
\end{table}

\begin{table}[h]
  \centering
      \begin{tabular}{c|cc|cc|ccc}
 & \multicolumn{2}{c}{GL1} & \multicolumn{2}{c}{GL2} & \multicolumn{2}{c}{GL3}\\
        $\Delta t$ &
        $T_{10^{-6}}/T_{10^{-5}}$ &
        $T_{10^{-7}}/T_{10^{-6}}$ &
                  $T_{10^{-6}}/T_{10^{-5}}$ &
          $T_{10^{-7}}/T_{10^{-6}}$ &
                    $T_{10^{-6}}/T_{10^{-5}}$ &
        $T_{10^{-7}}/T_{10^{-6}}$ \\
\hline
7200 & 1.006 & 1.243 & 1.006 & 1.393 &  1.031 & 1.208 \\
3600 & 1.465 & 1.159 & 1.014 & 1.361 &  1.007 & 1.009 \\
2400 & 1.252 & 1.183 & 1.023 & 1.188 &  1.266 & 1.148 \\
1200 & 1.292 & 1.147 & 1.091 & 1.316 &  1.375 & 1.160 \\
600 & 1.153 & 1.017 & 1.230 & 1.052 &   1.101 & 1.067 \\
300 & 1.221 & 1.078 & 0.957 & 1.230 &   1.217 & 1.093 \\

      \end{tabular}
      \caption{\label{tab:GLntol}Wallclock time ratios for
        Gauss-Legendre IRKs (stages 1,2,3) for various relative
        tolerances. Here, $T_{\tau}$ indicates the time to solution
        when nonlinear relative tolerance $\tau$ is used.
        }
  \end{table}

\begin{figure}
  \begin{center}
    \includegraphics[scale=0.4]{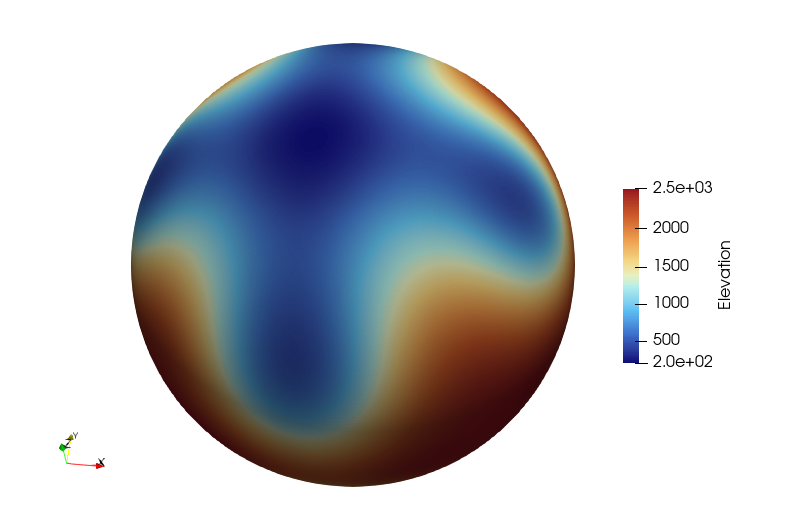}
  \end{center}
  \caption{\label{fig:w6_28}The elevation after 28 days in the Rossby-Haurwitz
    test case using $\Delta t=3600$ and the 2-stage Gauss-Legendre scheme.}
\end{figure}

\begin{figure}
  \begin{center}
    \includegraphics[scale=0.6]{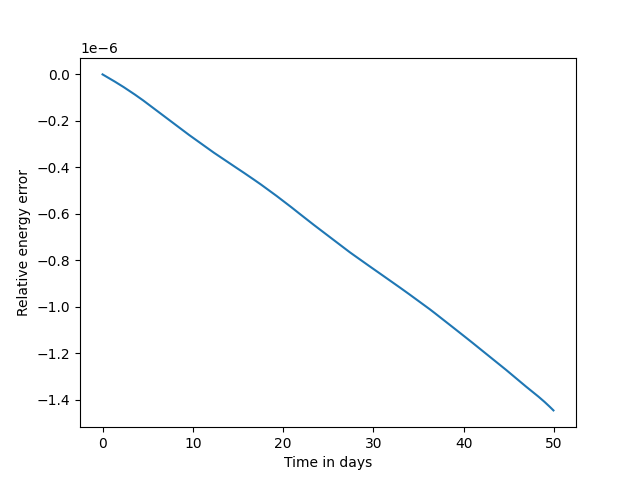}
    \includegraphics[scale=0.6]{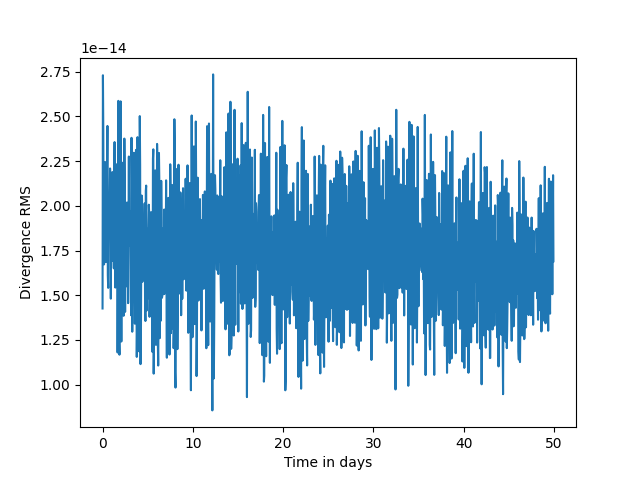}
  \end{center}
    \caption{\label{fig:w6energy}Plots of energy and the RMS of the
      divergence versus timestep over 50 days in the Rossby-Haurwitz
      test case using $\Delta t=3600$ and the 2-stage Gauss-Legendre scheme.}
\end{figure}

\begin{figure}
  \begin{center}
    \includegraphics[scale=0.4]{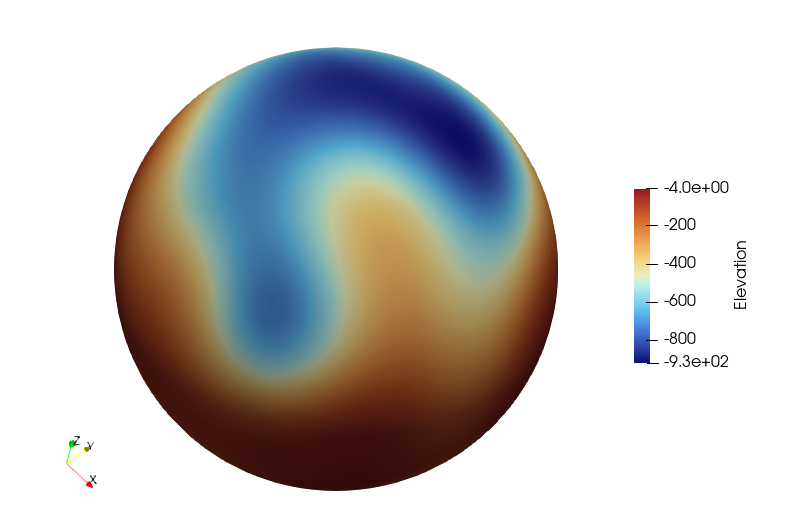}
    \includegraphics[scale=0.4]{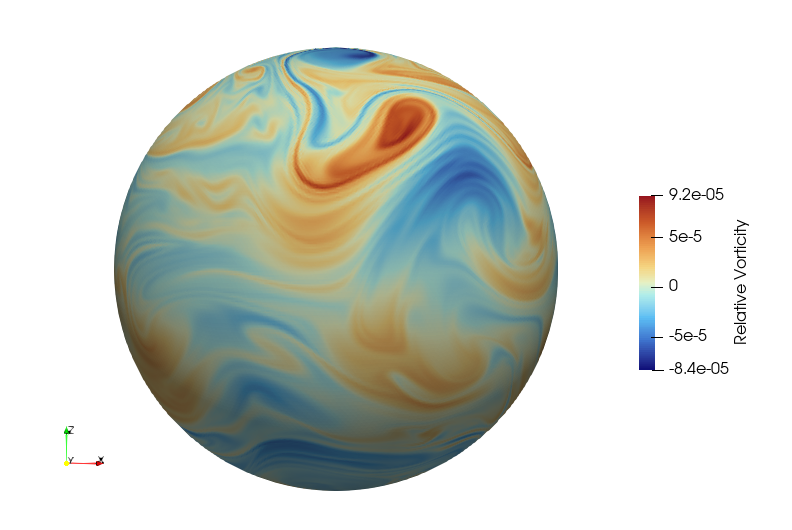}
  \end{center}
  \caption{\label{fig:w5_long}Top: The elevation after 15 days in the
    mountain test case using $\Delta t=3600$ and the 2-stage
    Gauss-Legendre scheme. Bottom: The vorticity for the same test
    after 50 days.}

\end{figure}

\begin{figure}
  \begin{center}
    \includegraphics[scale=0.6]{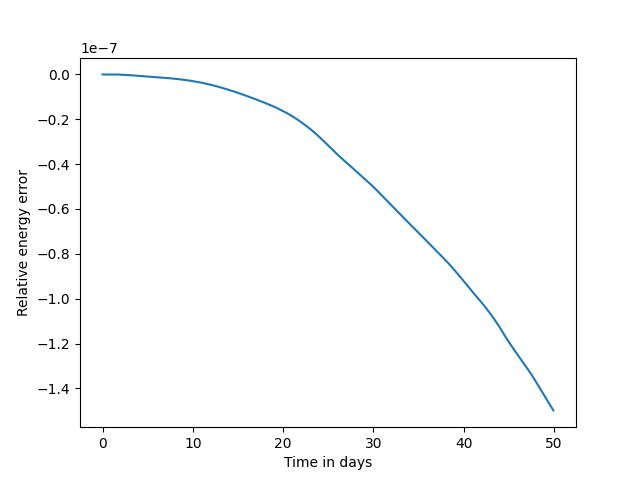}
    \includegraphics[scale=0.6]{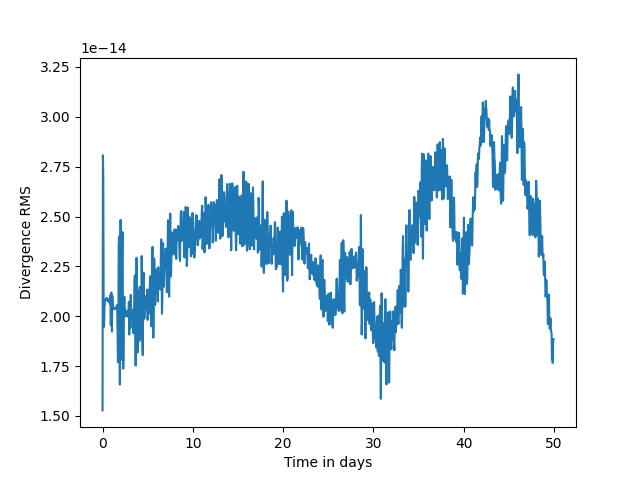}
  \end{center}
      \caption{\label{fig:w5energy}Plots of energy and the RMS of the
      divergence versus timestep over 50 days in the mountain
      test case using $\Delta t=3600$ and the 2-stage Gauss-Legendre scheme.}
\end{figure}

Finally, to demonstrate robustness of the numerical schemes we ran two
testcases for longer time periods, the mountain test case ``5'' in
addition to the Rossby-Haurwitz wave test case ``6'', both of
\citet{WILLIAMSON1992211}. We used the 2 stage Gauss-Legendre scheme for this,
with $\Delta t=3600$. Plots of the solution at longer times as well as
time series of energy and root mean square (RMS) of the divergence are shown in
Figures \ref{fig:w6_28}-\ref{fig:w5energy}. We observe that they are
both under control, indicating that the scheme is stable and not generating
spurious inertia-gravity waves.

 \section{Summary and outlook}
 \label{sec:summary}

An iterative solver strategy using multigrid and additive
  Schwarz smoothers applied to the monolithic coupled system for all
  variables at all stages was introduced for implicit Runge-Kutta
  methods applied to the rotating shallow water equations on the
  sphere, in conjunction with a compatible finite element
  discretisation. These patch systems are dense, and become reasonably
  large as $s$ increases ($144\times 144$ for $s=3$), which might be
  useful for modern computing architectures where floating point
  operations can be less important than data movement, so these
  methods are worth looking into now (although we do not explore this
  aspect here).

  We compared wallclock times and timestepping error
  for Gauss-Legendre and Radau IIA implicit collocation schemes with
  1, 2 and 3 stages with the ARK2 IMEX scheme. Our results demonstrate
  that these schemes achieve a comparable time to solution
    to the IMEX scheme for comparable accuracy, and can provide faster
  wallclock times if less accuracy is needed. Additionally, our solver
  approach could definitely benefit from further performance tuning,
  since the Firedrake multigrid implementation and patch
  frameworks are probably not optimal, whilst the direct solver used
  for the implicit solve in the IMEX scheme probably is,
    for problems of this size. Our results also demonstrate the
  capability of the combination of Irksome, Firedrake and PETSc to
  apply sophisticated solver strategies for IRKs.

These results give us confidence to go forward to examine IRK schemes
for three dimensional models. In that setting, direct solvers will not
be optimal for the IMEX schemes, and vertical line smoothers are
already required.
Hence, it is possible that IRKs will continue to be
  a competitive and viable option for this larger scale setting.
We will investigate this in future work.

\section*{Acknowledgements}
CC and JHC acknowledge the support of the Met Office via the Strategic Priorities Fund ExCALIBUR grant ``EX20-8 Exposing Parallelism: Parallel-in-Time''.

\section*{Code availability}
The version of the source code for this paper is archived at
\citet{colin_j_cotter_2026_20817145}. The version of Firedrake, PETSc
and other components is archived at
\citet{firedrake_zenodo_2026_20814010}.

\end{document}